\documentclass[reqno]{amsart}
\usepackage{latexsym}
\theoremstyle{plain}
\usepackage{txfonts}
\usepackage{enumerate}
\usepackage[sc]{mathpazo}
\def\XXint#1#2#3{{\setbox0=\hbox{$#1{#2#3}{\int}$}
     \vcenter{\hbox{$#2#3$}}\kern-.5\wd0}}

\newcommand{\lt}{\left}
\newcommand{\rt}{\right}
\newcommand{\nl}{\newline}

\newcommand{\nn}{\nonumber}

\newcommand{\lm}{\lambda}

\newcommand{\ep}{\epsilon}

\newcommand{\LI}{\mathcal{L}}

\newcommand{\ti}{\tilde}

\newcommand{\R}{\mathrm {I\!R}}

\newcommand{\dia}{\diamondsuit}

\newtheorem{a1}{Lemma}
\newtheorem{a2}{Theorem}

\theoremstyle{remark}

\begin{document}
\title[Differential inclusions, non-absolutely convergent integrals]{Differential inclusions, non-absolutely convergent integrals and the first theorem of complex analysis}
\author{Andrew Lorent}
\address{Mathematics Department\\University of Cincinnati\\2600 Clifton Ave.\\ Cincinnati OH 45221}
\email{lorentaw@uc.edu}
\maketitle

\begin{abstract} 

In the theory of complex valued functions of a complex variable arguably the first striking theorem is that pointwise differentiability 
implies $C^{\infty}$ regularity. As mentioned in Ahlfors \cite{ahl} there have been a number of studies \cite{porcon}, \cite{plunk1} proving this theorem without use of complex integration but at the 
cost of considerably more complexity.  In this note we will use the theory of non-absolutely convergent integrals to firstly give a 
very short proof of this result without complex integration and secondly (in combination with some elements of the theory of elliptic regularity) provide a far reaching generalization.
\end{abstract}

One of the first and most striking theorems about the analysis of complex valued functions of a complex variable is that merely from considering 
the class of pointwise complex differentiable functions on an open set we instantly find ourself in the category of $C^{\infty}$ functions.
\begin{a2}
\label{T0}
Given open set $\Omega\subset \mathbb{C}$. Suppose $f:\Omega\rightarrow \mathbb{C}$ is complex differentiable at every point. Then $f$ is $C^{\infty}$ on $\Omega$. 
\end{a2}

Typically Theorem \ref{T0} is proved via the method of complex integration. The first step is to prove that the integral of a differentiable function over the boundary of a rectangle inside a 
ball is zero, this was first proved by Goursat \cite{gou}. The existence of an anti-derivative is then concluded, Cauchy's integral formula follows and it is shown that you can differentiate through the integral of Cauchy's integral formula infinitely many times and hence 
the function is $C^{\infty}$. 

On the first paragraph of page 101 of Ahlfors's standard text \cite{ahl}, he writes that many important properties of analytic functions are 
difficult to prove without use of complex integration. Ahlfors states that only recently\footnote{\cite{ahl} was published in 1978} it has been possible to prove continuity of the gradient (or the existence of higher gradients) without the use of complex integration. He 
refers to articles of Plunkett \cite{plunk1}, and Porcelli and Connell \cite{porcon} both of which rely on a topological theorem of 
Whyburn \cite{why}. Ahlfors notes that both 
these proofs are much more complicated than the original proof. 

It appears the question of how to prove Theorem \ref{T0} without complex integration is a natural one and it was also asked by Luzin  \cite{landis} to two of his last students G.M.\ Adel'son-Vel'ski\v{i} and A.S.\ Kronrod. These two mathematicians produced a very interesting paper that solves the problem for domains bounded by a Jordan curve \cite{kon}. Amongst 
other ideas their method of proof anticipates what would now be called "bootstrapping" in the theory of elliptic PDE.    
Despite the fact the authors of this paper later became quite prominent scientists \footnote{ Indeed 
A.S. Kronod was in addition a decorated veteran of the Soviet counter offensive in the defence of Moscow in the 
winter of 1941. It appears he only returned to the study of mathematics after multiple wounds left him 
unable to fight \cite{landis}.} their paper has received very little attention, indeed according to \em mathscinet \rm this 
article will be the first one to cite it.

As mentioned, if we permit complex integration then Theorem \ref{T0} follows from Cauchy's integral formula which itself  
follows from Goursat's theorem. Generalizations of Goursat's theorem have a long history and one line of generalization provides an alternative proof 
of the theorem that in essence does not require anything from the theory of complex integration except its definition which is needed in the last step. This line of research was started by Montel \cite{mont} and 
further developed by Looman \cite{loo23} and Menchoff \cite{mench} (their theorem receives a very clear exposition in Saks \cite{saks}) and later 
by Tolstov \cite{tol}. Although not explicitly stated the method of proof of \cite{tol} was essentially to construct a Denjoy type integral 
to integrate the divergence of a differentiable vector field \cite{preiss6}. The explicit application of the theory of non-absolutely convergent integrals appears to have first been made by Mawhin \cite{mawhin1} who proved a general divergence theorem and used it to deduce Goursat's theorem, it is noted in \cite{mawhin1} that the fact Goursat's theorem is a corollary to  Green's theorem was known already to Bocher \cite{boch1} and de la Vall\'{e}e-Poussin \cite{val}. Later Jurkat and Nonnenmacher \cite{jur} used a kind of Perron integral in the plane to prove a generalization of Goursat's theorem due to Besicovitch \cite{besc}. Again they used similar strategy of first proving a general Green's theorem for class of vector fields that includes all differentiable vector fields, then deducing 
Gousat's theorem.

The purpose of this note is firstly to use the theory of non-absolutely convergent integrals to provide the shortest proof of Theorem 
\ref{T0}, the proof we provide is also independent of the theory of complex line integrals. Secondly by rephrasing Theorem \ref{T0} 
in terms of differential inclusions we will provide a far reaching generalization of this result by again applying the theory of 
non absolutely convergent integrals and some simple elliptic regularity estimates.

First some background. Note that the statement that $f(x+iy)=u(x,y)+iv(x,y)$ is pointwise complex differentiable on $\Omega\subset \mathbb{C}$ is equivalent to the 
statement that the vector valued function $\ti{f}(x,y)=(u(x,y),v(x,y))$ is pointwise differentiable and satisfies the differential 
inclusion 
\begin{equation}
\label{eq5}
D\ti{f}(x,y)\in \lt\{\lt(\begin{smallmatrix} a & -b \\ b & a\end{smallmatrix}\rt):a,b\in \R\rt\}=:\LI\text{ for any }x+iy\in \Omega. 
\end{equation}

The set $\LI$ has \em no rank-$1$ connections, \rm  by this we mean that if $A,B\in \LI$ and $\mathrm{rank}(A-B)=1$ then $A=B$. It turns out this is a 
crucial property that implies regularity of differential inclusions. We will establish the following generalization of Theorem 
\ref{T0}. 

\begin{a2}
\label{T1}
Suppose $u:\Omega\rightarrow \R^m$ is differentiable on open set $\Omega\subset \R^n$ and $L\subset M^{m\times n}$ is a 
subspace without rank-$1$ connections. If $D u(x)\in L$ for every $x\in \Omega$ then $u$ is real analytic in $\Omega$. 
\end{a2}

As noted one of the main ideas we will need to establish Theorem \ref{T1} is the use of the theory of non-absolutely continuous integrals. The 
main point about these integrals is that they allow us to integrate all derivatives of functions and thus they provide a stronger 
form of the fundamental theorem of calculus. On the real line this was accomplished by Denjoy \cite{den1}, \cite{den2} and later by 
Perron \cite{perron} by somewhat different ideas. Later in fifties in Kursweil \cite{kurz} and Henstock \cite{henstock} developed a different approach (sometimes known as the Gauge integral) that allowed for an easier generalization to higher dimensions. Recently there has been a strong line of generalization by Maly \cite{maly} who developed a theory of integration with respect to a distribution that he refers to as the CU integral. For our purposes 
we require any integral that in addition to basic linearity and finiteness properties satisfies 
\begin{equation}
\label{intpr1}
  \int \phi_{,i}(z) d\lm_n=0\text{ for any }\phi\in C_c(\R^n:\R), i\in \lt\{1,2,\dots n\rt\}.
\end{equation}

These properties are satisfied on the real line by any of the integrals mentioned. They are also satisfied in the plane by a version of 
the Perron integral \cite{jur}. They are satisfied in all dimension by the CU integral of \cite{maly}. The easiest integral to reference that satisfies these properties in all dimensions is the Gauge integral. In particular the divergence theorems of Mawhin \cite{mawhin1} and Pfeffer \cite{pf1} show that such an integral satisfies (\ref{intpr1}). For this reason we carefully detail the results on the Gauge integral that we need. In order to highlight the role of non-absolutely convergent integrals are playing 
in the arguments we will denote these integrals as $\int \dots d\lm_n$. Where as when the 
standard Lebesgue integral suffices we will denote these integrals by $\int \dots dz$.  \nl

\bf Acknowledgements. \rm Firstly I would like to thank David Preiss for suggesting the use of non-absolutely convergent integrals as a tool 
to further simplify my original simple proof of Theorem \ref{T0} and for many useful conversations about them. Secondly I would like to 
thank the participants of the UC geometric analysis seminar where I presented preliminary results including an approach to the conjecture 
that is now Theorem \ref{T1}. Examples provided by Robbie Buckingham, Gary Weiss and Catalin Dragan were partly responsible for the choice 
of a different approach to prove Theorem \ref{T1}. Finally I would like to thank the referee of an earlier version of 
this paper for suggesting a simplification to Lemma \ref{LL1}. I am also very grateful to the second referee who 
provided the reference \cite{kon} and pointed out the statement of the Theorem 2 can be strengthen to real analyticity by 
use of classic results detailed in \cite{morrey2}.

\subsection{Preliminaries} As in \cite{pf2} page 134 we define a cell 
in $\R^n$ to be a set of the form $A=\Pi_{i=1}^m A_i$ where $A_i:=\lt[\xi_i,\eta_i\rt]$ where we insist that the interior of this set (denoted $A^0$) is non empty. 

As outlined on page 210 \cite{pf2}, a figure is a finite (possibly empty) union of cells. On page 215 \cite{pf2} a set $T\subset \R^m$ is said to be 
thin if and only if it is the union of countably many sets whose $(m-1)$ dimensional Hausdorff measure is finite. We let $R_{*}(A,\lm_m)$ denote the set of Gauge integrable 
functions and $\bar{R_{*}}(A,\lm_m)$ be the set of extended real valued functions 
that agree with some function in $R_{*}(A,\lm_m)$ a.e.. The only 
result on the Gauge integral we need is the following powerful theorem (Theorem 11.4.10 \cite{pf2})  

\begin{a2} 
\label{TT1}
Let $T$ be a thin set and let $v$ be a continuous vector field on a figure $A$ that is differentiable in $A^{0}\backslash T$. Then 
$\mathrm{div} v\in R_{*}(A,\lm_m)$
\begin{equation}
\label{eqz5}
\int_{A} \mathrm{div} v \;d\lm_m=\int_{\partial A} v\cdot n \;d H^{m-1}.
\end{equation}
\end{a2}

Since $A$ is continuous on the figure $A$ which is a closed set, the right-hand side of (\ref{eqz5}) can be considered as any kind of surface integral, in particular as an integral with respect to Hausdorff measure. The left-hand side is a Gauge integral however and under the hypothesis has to be consider in that way. It is immediate from Theorem \ref{TT1} that the Gauge integral satisfies (\ref{intpr1}). 

%
%

\section{Elliptic estimates for differential inclusions}

As far as we are aware the following estimates are folk law, we learned of them from \cite{muller}, however the proofs are only sketched in \cite{muller} 
and the precise estimates we need are not stated so we prove the results in detail.  As noted in \cite{muller} many well known results in 
Elliptic regularity follow from the rigidity implied by the differential inclusion $Du\in L$ where $L$ is a subspace containing 
no rank-$1$ connections and $u$ is a Sobolev function. The reason our Theorem \ref{T1} is distinct from Theorem \ref{T3} (which will be 
stated later in this section) is that the 
space of differentiable functions is distinct from the space of Sobolev functions\footnote{Maly \cite{maly} defines a class of weakly differentiable functions that include both differentiable functions and Sobolev functions, this yields the intriguing possibility of an broader analysis of functions spaces.}, however in our opinion the heart of the matter is the powerful estimates of this section. Rather than give a proof based on the Fourier transform (specifically the fact that the 
Fourier transform of a gradient matrix is a rank-$1$ matrix) as sketched in \cite{muller}, we will use the Caccioppoli inequality for constant coefficient Elliptic systems \footnote{This was suggested by the referee of an earlier version of this paper.}. 

\begin{a1} 
\label{LL1}
Suppose $\Omega$ is a Lipschitz domain and $L$ is a subspace without rank-$1$ connections. If 
$v\in W^{1,2}(\Omega:\R^m)$ is such that $Dv\in L$ a.e.\ then for any sub-domain $U\subset \subset \Omega$ there 
exists constant $c=c(U)$ such that 
\begin{equation}
\label{eqzz1}
\int_{U} \lt|Dv\rt|^2 dz\leq c\int_{\Omega} \lt|v\rt|^2 dz.
\end{equation}
\end{a1}
\em Proof of Lemma \ref{LL1}. \rm 

\em Step 1. \rm We will show that if $v\in W^{1,2}(\Omega:\R^m)$ is such that  $Dv\in L$ a.e.\ then 
$v$ weakly satisfies a constant coefficient system 
\begin{equation}
\label{frvv1}
\sum_{\alpha=1}^n   \frac{\partial}{\partial x_{\alpha}} \lt(\sum_{j=1}^m\sum_{\beta=1}^n  \frac{\partial v_j}{\partial x_{\beta}}  A^{\alpha \beta}_{ij} \rt)=0
\end{equation}
where the  coefficients $A^{\alpha \beta}_{ij}$ are strongly elliptic in the sense that there 
exists a constant $\mu>0$ such that 
\begin{equation}
\label{frvv2}
\sum_{i=1}^m \sum_{\alpha=1}^n \sum_{j=1}^m\sum_{\beta=1}^n\ a_j b_{\beta} a_i b_{\alpha} A_{ij}^{\alpha\beta}
\geq \mu \lt|a\rt|^2\lt|b\rt|^2\text{ for every }a\in \R^m, b\in \R^n. 
\end{equation}

\em Proof of Step 1, \rm Given $a\in \R^m$, $b\in \R^n$ let $a\otimes b\in M^{m\times n}$ be the matrix 
whose $(i,j)$-th entry is $a_i b_j$. 
Let $A:M^{m\times n}\rightarrow M^{m\times n}$  be the orthogonal projection onto $L^{\perp}$. So if 
$A(a\otimes b)=0$ then $a\otimes b\in L$ which is a contradiction because $0\in L$ and so we would have a rank-$1$ connection in $L$. 
Now by homogeneity there exists constant $\lm>0$ such that 
\begin{equation}
\label{eqzz80}
\lt|A(a\otimes b)\rt|\geq \lm \lt|a\otimes b\rt|\text{ for every }a\in \R^m, b\in \R^n. 
\end{equation}
Now we can decompose $a\otimes b=A(a\otimes b)+P_{L}(a\otimes b)$ where $P_L$ is the orthognal projection onto 
$L$. So
\begin{eqnarray}
\label{uzua2}
a\otimes b: A(a\otimes b)&=& \lt(A(a\otimes b)+P_L(a\otimes b) \rt):A(a\otimes b)\nn\\
&=&\lt|A(a\otimes b)\rt|^2\overset{(\ref{eqzz80})}{\geq} \lm^2 \lt|a\otimes b\rt|^2\nn\\
&=&\lm^2 \lt|a\rt|^2 \lt|b\rt|^2\text{ for every }a\in \R^m, b\in \R^n.
\end{eqnarray}
Thus $A$ satisfies the strict Legrendre-Hadamard condition. 

Using notation consistent with \cite{gia} we will represent $A$ by a matrix in the following way. 
Let $e^i_{\alpha}$ be the matrix whose $(i,\alpha)$ entry is $1$ and $0$ everywhere else, then 
\begin{equation}
\label{uzua1}
A(e^j_{\beta})=\sum_{i=1}^{m}\sum_{\alpha=1}^n A_{ij}^{\alpha \beta} e_{\alpha}^i. 
\end{equation}
Writing (\ref{uzua2}) out in coordinates we obtain 
\begin{eqnarray}
\sum_{i=1}^m \sum_{\alpha=1}^n \sum_{j=1}^m\sum_{\beta=1}^n\ a_j b_{\beta} a_i b_{\alpha} A_{ij}^{\alpha\beta}
&=& a\otimes b:\lt(\sum_{i=1}^m\sum_{\alpha=1}^n \lt(\sum_{j=1}^m\sum_{\beta=1}^n\ a_j b_{\beta} A_{ij}^{\alpha\beta}\rt) e^i_{\alpha}\rt)\nn\\
&\overset{(\ref{uzua1})}{=}& a\otimes b:\lt(\sum_{j=1}^m \sum_{\beta=1}^n A(a_j b_{\beta} e^j_{\beta})\rt)\nn\\
&=& a\otimes b: A(a\otimes b) \overset{(\ref{uzua2})}{\geq} \lambda^2 \lt|a\rt|^2\lt|b\rt|^2\text{ for every }a\in \R^m, b\in \R^n. \nn
\end{eqnarray}
So taking $\mu=\lm^2$ inequality (\ref{frvv2}) is satisfied.  

Now note $A(Dv)=0$ a.e.\ and so $v$ weakly satisfies $\mathrm{div}(A(Dv))=0$. And as 
\begin{eqnarray}
A(Dv)&=&\sum_{j=1}^m\sum_{\beta=1}^n \frac{\partial v_j}{\partial x_{\beta}} A(e^j_{\beta})\overset{(\ref{uzua1})}{=}\sum_{i=1}^m\sum_{\alpha=1}^n \lt( \sum_{j=1}^m\sum_{\beta=1}^n \frac{\partial v_j}{\partial x_{\beta}} A^{\alpha \beta}_{ij}  \rt) e^i_{\alpha}.\nn 
\end{eqnarray}
Now letting $\lt[M\rt]_i$ denote the $i$-th row of matrix so 
$$
\mathrm{div}\lt(\lt[A(Dv)\rt]_i\rt)=\sum_{\alpha=1}^n   \frac{\partial}{\partial x_{\alpha}} \lt(\sum_{j=1}^m\sum_{\beta=1}^n  \frac{\partial v_j}{\partial x_{\beta}}  A^{\alpha \beta}_{ij} \rt)=0
$$
thus (\ref{frvv1}) is satisfied. \nl

\em Step 2. \rm We will show $v$ satisfies (\ref{eqzz1}). 

\em Proof of Step 2. \rm Recall the $M^{nm\times nm}$ matrix of coefficients of the orthogonal projection $A$ are denoted $A^{\alpha \beta}_{ij}$.  
By (\ref{frvv2}) the matrix $\lt(A^{\alpha \beta}_{ij}\rt)$ satisfies the matrix form of the strict Legrendre-Hadamard condition given by (2.2) p76 \cite{gia} and we have that $\mathrm{div}\lt(\lt[A(Dv)\rt]\rt)$ is the system of $m$ equations given by (2.1) p76 \cite{gia}. Hence 
by the Caccioppoli inequality given by Proposition 2.1. \cite{gia} function $v$ satisfies (\ref{eqzz1}). $\Box$


\begin{a2}
\label{T3}
Suppose $u\in W^{1,2}(\Omega:\R^m)$ with $Du(x)\in L$ a.e.\ $x\in \Omega$ then $u$ is real analytic in $\Omega$.
\end{a2}

\em Proof of Theorem \ref{T3}. \rm 

\em Step 1. \rm We will show $u$ is $C^{\infty}$. 

\em Proof of Step 1. \rm Let $i\in \lt\{1,2,\dots n\rt\}$. Let 
$\delta>0$ and $\Pi_0:=\Omega\backslash N_{\delta}(\Omega)$, $\Pi_1:=\Pi_0\backslash N_{\delta}(\Pi_0)$. For any $h\in (0,\delta)$ we 
have $w_h(x):=\frac{u(x+he_i)-u(x)}{h}$ we have $D w_h(x)\in L$ for a.e.\ $x\in \Pi_0$. So by Lemma \ref{LL1} we have 
$$
\int_{\Pi_1} \lt|D w_h \rt|^2 dz\leq c\int_{\Pi_0} \lt|w_h\rt|^2 dz.
$$
Since $u\in W^{1,2}(\Omega:\R^n)$, $w_h\overset{L^2(\Pi_0)}{\rightarrow} \frac{\partial u}{\partial x_i}$ and 
$\frac{\partial u}{\partial x_i}\in L^2(\Omega)$ so 
$\int_{\Pi_0} \lt|w_h\rt| dz\leq c$ for all small enough $h>0$. Since $i$ is arbitrary this implies $Du\in W^{1,2}(\Pi_1)$ and 
$$
\int_{\Pi_1} \lt|D^2 u\rt|^2 dz\leq c \int_{\Pi_0} \lt|D u\rt|^2 dz.
$$
Now as $D w_h(z)\in L$ for a.e.\ $z$ and $D w_h\overset{L^2(\Pi_0)}{\rightarrow} D u_{,i}$ so we know $Du_{,i}\in L$ for 
a.e.\ $x\in \Pi_0$. So we can repeat the argument for the function $u_{,i}$ and gain control of the third order derivatives.  By repeating 
in this way the result follows from Sobolev embedding theorem. \nl

\em Step 2. \rm We will show $u$ is real analytic.

\em Proof of Step 2. \rm By Step 1 of Lemma \ref{LL1} $u$ weakly satisfies the constant coefficient system 
(\ref{frvv1}). Since by Step 1 of this Lemma $u$ is $C^{\infty}$ so it strongly satisfies this system. Now from Definition 
6.5.1 and equations (6.5.10), (6.6.4) of \cite{morrey2}, by (\ref{frvv2}) this system is "strongly elliptic" system in the sense of \cite{morrey2} 
and it is a system whose "operators $L_{jk}$" have constant constant coefficients and hence are analytic. Thus by Theorem 6.6.1 of \cite{morrey2} $u$ is real analytic in $\Omega$. $\Box$

\section{Proof of Theorem \ref{T1}.}

\em Step 1. \rm We will show that denoting $u^{\ep}:=u*\rho_{\ep}$ we have 
\begin{equation}
\label{eqzz52}
D u^{\ep}(x)\in L\text{ for any }x\in \Omega\backslash N_{\ep}(\partial \Omega). 
\end{equation}

\em Proof of Step 1. \rm Given $\psi\in C^{\infty}(\Omega:M^{m\times n})$. Define the $m$ vector 
\begin{equation}
\label{eqzz30}
\lt[\Lambda_z(\psi)\rt](z):=\lt(\begin{matrix} \sum_{j=1}^n \frac{\partial \psi_{1,j}(z)}{\partial z_j}\nn\\
 \sum_{j=1}^n \frac{\partial \psi_{2,j}(z)}{\partial z_j} \nn\\
\dots \nn\\
\sum_{j=1}^n \frac{\partial \psi_{m,j}(z)}{\partial z_j} \nn
 \end{matrix}\rt).
\end{equation}
Let $w\in C^1(\Omega':\R^m)$ where $\Omega'\subset \Omega$ and 
$\mathrm{spt}\psi\subset \Omega'$. Note that for any $i\in \lt\{1,2,\dots m\rt\}$
\begin{eqnarray}
\label{eqzz201}
\mathrm{div}\lt(w_i(z) \psi_{i,1}(z), w_i(z) \psi_{i,2}(z),\dots  w_i(z) \psi_{i,n}(z)\rt)=\sum_{j=1}^n \lt(w_{i,j}(z)\psi_{i,j}(z)+
w_{i}(z)\frac{\partial \psi_{i,j}}{\partial z_j}(z)\rt).
\end{eqnarray}
Note 
\begin{equation}
\label{eqzz21}
\int \mathrm{div}\lt(w_i(z) \psi_{i,1}(z), w_i(z) \psi_{i,2}(z),\dots  w_i(z) \psi_{i,n}(z)\rt) dz=0. 
\end{equation}
So putting this together with (\ref{eqzz201}) we have 
\begin{eqnarray}
\label{eqzz20}
\int  \sum_{i=1}^m\sum_{j=1}^n w_{i,j}(z)\psi_{i,j}(z) dz&=&
-\int \sum_{i=1}^m\sum_{j=1}^n w_{i}(z)\frac{\partial \psi_{i,j}}{\partial z_j}(z) dz\nn\\
&=&-\int w(z)\cdot \lt(\lt[\Lambda_z (\psi)\rt](z)\rt) dz.
\end{eqnarray}

Let $\phi\in C^{\infty}_c(\Omega:M^{m\times n})$ with $\phi(z)\in L^{\perp}$ for all $z\in \Omega$. Take $\ep>0$ sufficiently small to that $\mathrm{spt} \phi\subset \Omega\backslash N_{\ep}(\partial \Omega)$ 
\begin{eqnarray}
\label{eqzz41.5}
\int u^{\ep}(z)\cdot \lt[\Lambda_z (\phi)\rt](z) dz&=&\int \lt(\int u(x)\rho_{\ep}(x-z) dx\rt) \cdot \lt(\lt[\Lambda_z (\phi)\rt](z)\rt) dz\nn\\
&=&\int \int u(x)\rho_{\ep}(x-z)\cdot  \lt(\lt[\Lambda_z (\phi)\rt](z)\rt) dz dx\nn\\
&=&\int \int \sum_{i=1}^m \sum_{j=1}^n \lt(u_i(x)\rho_{\ep}(x-z)\rt)\frac{\partial \phi_{i,j}(z)}{\partial z_j} dz dx
\end{eqnarray}
But notice 
\begin{eqnarray}
\label{ffeq9}
\int \sum_{i=1}^m \sum_{j=1}^n \lt(u_i(x)\rho_{\ep}(x-z)\rt)\frac{\partial \phi_{i,j}(z)}{\partial z_j} dz&=&
\sum_{i=1}^m \sum_{j=1}^n \int u_i(x) \frac{\partial \rho_{\ep}}{\partial z_j}(x-z)  \phi_{i,j}(z) dz\nn\\
&=&\sum_{i=1}^m \sum_{j=1}^n -\int u_i(x) \frac{\partial \rho_{\ep}}{\partial x_j}(x-z) \phi_{i,j}(z) dz.
\end{eqnarray}
Now note by (\ref{intpr1}) $\int \frac{\partial }{\partial x_j}\lt(u_i(x)\rho_{\ep}(x-z)\rt) d\lm_n(x)=0$. So 
\begin{equation}
\label{ffeq11}
\int \frac{\partial u_i}{\partial x_j}(x)\rho_{\ep}(x-z) d\lm_n(x)=-\int u_i(x) 
\frac{\partial \rho_{\ep}}{\partial x_j}(x-z) dx.
\end{equation}
Putting these things together we have 
\begin{eqnarray}
\label{eqzz41}
\int u^{\ep}(z)\cdot \lt[\Lambda_z (\phi)\rt](z) dz&\overset{(\ref{ffeq9}),(\ref{eqzz41.5})}{=}&-\int \int \sum_{i=1}^m \sum_{j=1}^n  
 u_i(x) \frac{\partial \rho_{\ep}}{\partial x_j}(x-z) \phi_{i,j}(z) dx dz\nn\\
&\overset{(\ref{ffeq11})}{=}&\int \int \sum_{i=1}^m \sum_{j=1}^n \frac{\partial u_i}{\partial x_j}(x) \rho_{\ep}(x-z) \phi_{i,j}(z) d\lm_n(x) dz\nn\\
&=&\int D u(x):\phi*\rho_{\ep}(x) d\lm_n(x)\nn\\
&=&0
\end{eqnarray}
since $\phi*\rho_{\ep}(x)\in L$ for any $x$.

Thus by (\ref{eqzz20}) equation (\ref{eqzz41}) this implies that $\int D u^{\ep}(z):\phi(z) dz=0$ for any $\phi\in C^{\infty}(\Omega:L^{\perp})$ which 
implies (\ref{eqzz52}).\nl

\em Proof of Theorem \ref{T1} completed. \rm Let $U\subset \subset V\subset \subset \Omega$. Let $\ep_{n}\rightarrow 0$. Now by Lemma \ref{LL1} and Step 1 we know 
\begin{equation}
\label{eqzz60}
\int_{U} \lt|D u^{\ep_n}\rt|^2 dz\leq c \int_{V} \lt|u^{\ep_n}\rt|^2 dz.
\end{equation}
Since $u^{\ep_n}\overset{L^2(V)}{\rightarrow} u$ so there exists constant $c_0$ such that $\int_{U} \lt|u^{\ep_n}\rt|^2 dx\leq c_0$ for 
all $n$. So $u^{\ep_n}$ is a bounded sequence in $W^{1,2}(U)$ and we can extract a weakly converging subsequence $u^{\ep_{k_n}}$ 
that converge to $u$. So $u\in W^{1,2}(U)$ and $Du\in L$ a.e.\ in $\Omega$. So by Theorem \ref{T3} this implies $u$ is analytic in $U$. As $U$ is an arbitrary subset 
of $\Omega$ this implies $u$ in analytic in $\Omega$. $\Box$

\section{A simple proof of Theorem \ref{T0}}

Theorem \ref{T0} is a very special case of Theorem \ref{T1} that we have just proved. However 
in this section we will show the theory of non-absolutely convergent integrals and Weyl's lemma allows us a very short direct proof of 
Theorem \ref{T0}.  

\begin{a2} 
\label{T2}
Given open set $\Omega\subset \mathbb{C}$. Suppose $f:\Omega\rightarrow \mathbb{C}$ is a complex differentiable at every point. 
Let $u,v:\Omega\rightarrow \R$ be defined by $f(x+iy)=u(x,y)+iv(x,y)$. We will show $u,v$ weakly satisfy Laplace's equation. 
\end{a2}
\em Proof of Theorem \ref{T2}. \rm Firstly since $f(x+iy)=u(x,y)+iv(x,y)$ is complex differentiable, 
for any vector $(r,s)$ we have 
$$
f'(x+iy)=\lim_{h\rightarrow 0} \frac{u\lt((x,y)+h(r,s)\rt)-u(x,y)}{h}
+i\frac{v\lt((x,y)+h(r,s)\rt)-v(x,y)}{h}
$$   
and so $u$, $v$ are pointwise differentiable in $\Omega$. Now definition of complex differentiability we have $f'(x+iy)=u_x(x,y)+iv_x(x,y)=-iu_y(x,y)+v_y(x,y)$. So $u,v$ satisfy the Cauchy Riemann equation. Take $\phi\in C_{c}^{\infty}(\Omega)$. Note $\mathrm{div}(-\phi_y v,\phi_x v)=-\phi_y v_x+\phi_x v_y$ and $(-\phi_y v,\phi_x v)$ is a differentiable vector field. 

Take $R>0$ such that $\Omega\subset Q_R(0)$, by (\ref{intpr1}) (or by Theorem \ref{TT1})  
we have that 
\begin{equation}
\label{eq1}
0=\int_{Q_R(0)} \mathrm{div}(-\phi_y v,\phi_x v)\; d\lm_2=\int_{\Omega} -\phi_y v_x+\phi_x v_y\; d\lm_2 =\int_{\Omega} \phi_y u_y+\phi_x u_x\; d\lm_2.
\end{equation}

Now note that vector field $(\phi_x u, \phi_y u)$ is differentiable and  
\begin{eqnarray}
\label{eq2}
0&=&\int_{Q_R(0)} \mathrm{div}(\phi_x u,\phi_y u)\; d\lm_2 =\int_{\Omega} \triangle \phi u+\phi_x u_x+\phi_y u_y \;d\lm_2 \overset{(\ref{eq1})}{=}\int_{\Omega} \triangle \phi u\; d\lm_2.
\end{eqnarray}
Note that the integral on the right-hand side of (\ref{eq2}) is the integral of a continuous function as an integral it is equal to the 
Lebesgue integral. So specifically we have shown 
\begin{equation}
\label{eq3}
\int_{\Omega} \triangle \phi u\; dz=0\text{ for any }\phi\in C_{c}^{\infty}(\Omega)
\end{equation}
So $u$ weakly satisfy Laplace's equation. 

Arguing in the same way. Take $\phi\in C_{c}^{\infty}(\Omega)$. 
Note the vector field $(\phi_y u,-\phi_x u)$ is differentiable so applying (\ref{intpr1}) (or Theorem \ref{TT1}) 
we have 
\begin{equation}
\label{eq4}
0=\int_{Q_R(0)} \mathrm{div}(\phi_y u,-\phi_x u)\; d\lm_2=\int_{\Omega} \phi_y v_y+\phi_x v_x\; d\lm_2.
\end{equation}
And again by differentiability of the vector field $(\phi_x v,\phi_y v)$ we have 
\begin{equation}
\label{eq5}
0=\int_{Q_R(0)} \mathrm{div}(\phi_x v,\phi_y v)\;d\lm_2=\int_{\Omega} \triangle \phi v+\phi_x v_x+\phi_y v_y\; d\lm_2
\overset{(\ref{eq4})}{=}\int_{\Omega} \triangle \phi v dz
\end{equation}
So again $v$ weakly satisfies Laplace's equation. $\Box$

\subsection{Proof of Theorem \ref{T0} continued.} We write $f(x+iy)=u(x,y)+iv(x,y)$, by Theorem \ref{T2} $u$ and $v$ weakly satisfy 
Laplace's equation. So by Weyl's lemma $u,v$ are $C^{\infty}$. So as $f'(x+iy)=u_x(x,y)+iv_x(x,y)$ and $u_x$, $v_x$ satisfy 
the Cauchy Riemann equations so $f'$ is complex differentiable. In the same way all orders of derivative of $f$ exists.  $\Box$

%
%

\section{Appendix} 
\subsection{Weyl's lemma} Weyl's lemma is well known but how easy and elementary its proof is perhaps is less well dissipated. For completeness 
we briefly outline the three main points.  

Firstly if we have a weakly harmonic function $u$ defined on $\Omega$ then the 
convolution $v=u*\rho_{\ep}$ is harmonic on $\Omega\backslash N_{\ep}(\partial \Omega)$, this follow from the definition of 
weak harmonicity and differentiating through the integral of the convolution. 

The second point is that it follows from the first point that weakly harmonic functions are uniformly approximated by 
harmonic functions and hence must satisfy the mean value theorem. 

The third point is that if function $u$ satisfies the mean value theorem then letting $\rho$ be a radial symmetric convolution 
kernel we have 
\begin{eqnarray}
u*\rho_{\ep}(z)&=&\int u(x)\rho_{\ep}(z-x) dx=\int_{0}^{\infty} \int_{\partial B_r(z)} u(x)\rho_{\ep}(z-x) dH^{n-1}x dr\nn\\
&=&\int_{0}^{\infty} \rho_{\ep}(z-re_1) u(z) H^{n-1}(\partial B_r(z)) dr =u(z)\int \rho_{\ep}(z-x) dx =u(z).\nn
\end{eqnarray}
So $u$ is $C^{\infty}$ and this completes the sketch.

\end{document}